\documentclass{amsart}
\usepackage{amssymb}
\usepackage{url}
\usepackage{algorithmic}

\newcommand\blfootnote[1]{%
  \begingroup
  \renewcommand\thefootnote{}\footnote{#1}%
  \addtocounter{footnote}{-1}%
  \endgroup
}

\newcommand{\re}{\mathrm{Re}}% partie reelle
\newtheorem{theorem}{Theorem}[section]

\theoremstyle{definition}

\theoremstyle{remark}

\numberwithin{equation}{section}

\begin{document}

\title{A New Effective Asymptotic Formula for the Stieltjes Constants}

\author[L. Fekih-Ahmed]{Lazhar Fekih-Ahmed}
\address{\'{E}cole Nationale d'Ing\'{e}nieurs de Tunis, BP 37,
Le Belv\'{e}d\`{e}re 1002 , Tunis, Tunisia}

\curraddr{} \email{lazhar.fekihahmed@enit.rnu.tn}
\thanks{}

%    \subjclass is required.


\date{June 30, 2014}

\dedicatory{}

%    Abstract is required.
\begin{abstract}
We derive a new integral formula for the Stieltjes constants. The
new formula permits easy computations using an effective
asymptotic formula. Both the sign oscillations and the leading
order of growth are provided. The formula can also be easily
extended to some generalized Euler constants.
\end{abstract}

\maketitle

\bibliographystyle{amsplain}

%    Section headings
\section{Introduction}\label{sec1}
The Stieltjes constants $\gamma_n$ are defined as  the
coefficients of Laurent series expansion of the Riemann zeta
function at $s=1$ ~\cite{briggs:zeta}:

\begin{equation}\label{sec1-eq1}
\zeta(s)=\frac{1}{s-1}+\sum_{n=0}^{\infty}\frac{(-1)^{n}}{n!}\gamma_{n}(s-1)^{n},
\end{equation}

where $\gamma_0=0.5772156649$ is known as Euler's constant.

Exact and asymptotic formulas as well as upper bounds for the
Stieltjes constants have been a subject of research for many
decades
\cite{briggs:zeta,briggs:chowla,coffey:stieltjes,johansson:hurwitz,knessl:coffey,kreminski:stieltjes,matsuoka:zeta,
mitrovic:stieltjes,zhang:williams}. The approach to estimate the
Stieltjes constants is always deterministic except the paper
\cite{adell:stieltjes} where a probabilistic approach is
undertaken. The main reason to estimate the Stieltjes constants is
that these constants and their generalization , known as
Generalized Euler constants, have many applications in number
theory.

This paper is a continuation of this line of research. We will
give a new effective  asymptotic formula for the Stieltjes
constants. With the formula we obtain the sign oscillations and
the leading order of growth of the Stieltjes constants. We will
show that our results match those of \cite{knessl:coffey} which
may be considered very accurate compared to other results.

\section{A New  Formula for The Stieltjes
Constants}\label{sec2}

Let $\phi(t)$ be the real function  defined by

\begin{equation}\label{sec2-eq1}
\phi(t)=\frac{d}{dt}\frac{-te^{-t}}{1-e^{-t}}=\frac{te^{t}}{(e^{t}-1)^2}-\frac{1}{e^{t}-1}.
\end{equation}

In a previous article we have obtained the following integral
representation of the Riemann zeta function:

\begin{theorem}[\cite{fekih:hurwitz}]\label{sec2-thm1}
With  $ \phi(t)$ as above, and for  all $s$ such that $\re(s)>-k$,
we have
\begin{equation}\label{sec2-eq3}
(s-1)\zeta(s)
=\frac{(-1)^k}{\Gamma(s+k)}\int_{0}^{\infty}\frac{d^{k}\phi(t)}{dt^k}
t^{s+k-1}\,dt.
\end{equation}
\end{theorem}

If we chose $k=1$ and we call

\begin{equation}\label{sec2-eq4}
\mu(t)=-\frac{d\phi}{dt}=\frac{d^2}{dt^2}\frac{te^{-t}}{1-e^{-t}}=-\frac{(2+t)e^{t}}{(e^{t}-1)^2}+\frac{2te^{2t}}{(e^{t}-1)^3},
\end{equation}

then Theorem~\ref{sec2-thm1} provides the following formula valid
for all $s$ such that $\re(s)>-1$

\begin{equation}\label{sec2-eq5}
s(s-1)\zeta(s)\Gamma(s) =\int_{0}^{\infty}\mu(t) t^{s}\,dt.
\end{equation}

If we now replace $s$ by $1-s$ in equation (\ref{sec2-eq5}) with
the assumption that $\re(1-s)<2$, we get

\begin{equation}\label{sec2-eq6}
s(s-1)\zeta(1-s)\Gamma(1-s) =\int_{0}^{\infty}\mu(t) t^{1-s}\,dt.
\end{equation}

The functional equation for the Riemann zeta function states that

\begin{equation}\label{sec2-eq7}
\zeta(s)=2(2\pi)^{s-1}\sin(\frac{\pi s}{2})\zeta(1-s)\Gamma(1-s).
\end{equation}

Multiplying both sides of the last equation by $s(s-1)$ and using
(\ref{sec2-eq6}), we obtain

\begin{eqnarray}\label{sec2-eq8}
s(s-1)\zeta(s)&=&2(2\pi)^{s-1}\sin(\frac{\pi
s}{2})s(s-1)\zeta(1-s)\Gamma(1-s)\nonumber\\
&=&2(2\pi)^{s-1}\sin(\frac{\pi s}{2})\int_{0}^{\infty}\mu(t)
t^{1-s}\,dt.
\end{eqnarray}

By observing that $(2\pi)^{s-1}=e^{(s-1)\log(2\pi)}$, that
$t^{1-s}=e^{-(s-1)\log(t)}$ and that $2\sin(\frac{\pi
s}{2})=2\cos(\pi\frac{(s-1)}{2})=e^{i\pi\frac{
(s-1)}{2}}+e^{-i\pi\frac{ (s-1)}{2}}$, we can rewrite
(\ref{sec2-eq8}) as

\begin{equation}\label{sec2-eq9}
s(s-1)\zeta(s)=\int_{0}^{\infty}\mu(t)
\big[e^{(s-1)(a-\log(t))}+e^{(s-1)(\bar{a}-\log(t))}\big]\,dt,
\end{equation}

where $a$ is the fixed complex number
$a=\log(2\pi)+i\frac{\pi}{2}$.

Finally, since the left hand side of (\ref{sec2-eq9}) is analytic
at $s=1$ it has a Taylor series expansion

\begin{equation}\label{sec2-eq10}
s(s-1)\zeta(s) =\sum_{n=0}^{\infty}\mu_{n}(s-1)^{n},
\end{equation}

where the coefficients $\mu_n$ are given by

\begin{eqnarray}\nonumber
\mu_n&=& \frac{1}{n!}\lim_{s\to 1}\frac{d^n}{d s^n}\big \{s(s-1)\zeta(s)\big \}\\
&=& \frac{1}{n!}\int_{0}^{\infty}\mu(t)\lim_{s\to
1}\frac{d^n}{d s^n} \big \{ e^{(s-1)(a-\log(t))}+e^{(s-1)(\bar{a}-\log(t))}\big \}\,dt\nonumber\\
&=&\frac{1}{n!}\int_{0}^{\infty}\mu(t) \big \{
(a-\log(t))^n+(\bar{a}-\log(t))^n\big \}\,dt.\label{sec2-eq11}
\end{eqnarray}

This gives our first main result\footnote{Note that the
coefficients $\mu_n$ and the integral formula of $s(s-1)\zeta(s)$
are as important as the Stieltjes constants and the function
$(s-1)\zeta(s)$. In fact, like the Riemann $\xi(s)=\tfrac{1}{2}
s(s-1)\pi^{-s/2}\Gamma\left(\tfrac{1}{2}s\right)\zeta(s)$
function, $s(s-1)\zeta(s)$ possess some symmetry and can play an
important role in the theory of the Riemann zeta function.}:

\begin{theorem}\label{sec2-thm2}
With  $ \mu(t)$ and the constant $a$ defined as above, the
coefficients $\mu_{n}$ are given by
\begin{equation}\label{sec2-eq12}
\mu_n=\frac{2}{n!}\int_{0}^{\infty}\mu(t)\re
\big\{(a-\log{t})^{n}\big\}\,dt.
\end{equation}

\end{theorem}

Once we have the coefficients $\mu_n$ of the power series for
$s(s-1)\zeta(s)$, the Stieltjes coefficients $\gamma_n$ can be
calculated using power series multiplication:

\begin{equation}\label{sec2-eq13}
(s-1)\zeta(s)=\sum_{n=0}^{\infty}(-1)^{n}(s-1)^{n}\times\sum_{n=0}^{\infty}\mu_n(s-1)^n
\end{equation}

since
\begin{equation}\label{sec2-eq14}
\frac{1}{s} =\sum_{n=0}^{\infty}(-1)^{n}(s-1)^{n}.
\end{equation}

This immediately yields

\begin{equation}\label{sec2-eq15}
\zeta(s)=\frac{1}{s-1}+\sum_{n=1}^{\infty}\left\{\sum_{k=0}^{n}(-1)^{n-k}\mu_k\right\}(s-1)^{n-1};
\end{equation}

therefore,

\begin{equation}\label{sec2-eq16}
\gamma_{n}=-n!\sum_{k=0}^{n+1}(-1)^{k}\mu_k=-
n!\int_{0}^{\infty}2\mu(t)\re
\big\{\sum_{k=0}^{n+1}\frac{(\log{t}-a)^{k}}{k!}\big\}\,dt.
\end{equation}

The last formula can be simplified even further. Indeed, the sum
inside the integral is a truncated sum of the exponential series
$e^{\log{t}-a}=te^{-a}$. This yields,

\begin{eqnarray}
\gamma_{n}&=&-n!\int_{0}^{\infty}2\mu(t)\re
\big\{\sum_{k=0}^{n+1}\frac{(\log{t}-a)^{k}}{k!}\big\}\,dt\nonumber\\
&=&-n!\int_{0}^{\infty}2\mu(t)\re \big\{
te^{-a}-\sum_{k=n+2}^{\infty}\frac{(\log{t}-a)^{k}}{k!}\big\}\,dt\nonumber\\
&=&n!\int_{0}^{\infty}2\mu(t)\re \big\{
\sum_{k=n+2}^{\infty}\frac{(\log{t}-a)^{k}}{k!}\big\}\,dt\label{sec2-eq17}
\end{eqnarray}

since $\re\{e^{-a}\}=\re\{\frac{-i}{2\pi}\}=0$. Hence, with

\begin{equation}\label{sec2-eq17bis}
I(n)=\int_{0}^{\infty}\mu(t) (\log{t}-a)^{n}\,dt,
\end{equation}

we can write

\begin{align}\label{sec2-eq17bb}
\gamma_{n}&=
n!\left[\frac{I(n+2)}{(n+2)!}+\frac{I(n+3)}{(n+3)!}+\ldots\right]\\
&=n!\left[(-1)^{n+2}\mu_{n+2}+(-1)^{n+3}\mu_{n+3}+\ldots\right],\label{sec2-eq17bbb}
\end{align}

and we have our second main result:

\begin{theorem}\label{sec2-thm3}
With  $\mu(t)$ and $\mu_n$ defined as above, the Stieltjes
constants are given by
\begin{equation}\label{sec2-eq18}
\gamma_{n}=
n!(-1)^n\left[\mu_{n+2}-\mu_{n+3}+\mu_{n+4}-\ldots\right].
\end{equation}
\end{theorem}

We do not know yet that the leading term
$n!(-1)^n\mu_{n+2}=\frac{n!}{(n+2)!}I(n+2)$ is the dominant term
for approximating $\gamma_{n}$.  All we know for now is that the
$\mu_n$'s are the Taylor coefficients of an entire function and
that $\mu_n\to 0$ as $n\to\infty$. In the next section, we will
show that  $|\mu_{n+2}|\gg |\mu_{n+3}|\gg \ldots$ for large $n$ so
that $\left\{|\mu_{k}|\right\}_{n+2}^{\infty}$ form an asymptotic
sequence\footnote{We write $f(n) \gg g(n)$, or $f$ is ``much
greater than" $g$, if $g = o(f)$ as $n\to \infty$.}. This implies
that $\gamma_{n}$ can be written as

\begin{equation}\label{sec2-eq18bis}
\gamma_{n}=\frac{1}{(n+1)(n+2)}\int_{0}^{\infty}2\mu(t)\re
\big\{(\log{t}-a)^{n+2}\big\}\,dt+\text{higher order terms},
\end{equation}

and that the leading  term  provides an asymptotic approximation
of $\gamma_{n}$:

\begin{equation}\label{sec2-eq19}
\gamma_{n}\approx \frac{1}{(n+1)(n+2)}\int_{0}^{\infty}2\mu(t)\re
\big\{(\log{t}-a)^{n+2}\big\}\,dt=n!(-1)^n\mu_{n+2}.
\end{equation}

Therefore, Theorem~\ref{sec2-thm3} permits an  asymptotic
expansion of the constants  $\gamma_{n}$. It is the subject of the
next section.

\section{Asymptotic Estimates of The Stieltjes Constants}\label{sec3}

This section is dedicated to approximating the complex-valued
integral

\begin{equation}\label{sec3-eq1}
I(n)=\int_{0}^{\infty}\mu(t) (\log{t}-a)^{n}\,dt.
\end{equation}

There are mainly two methods used for the asymptotic evaluation of
complex integrals of the form (\ref{sec3-eq1}) when $n$ is large:
the steepest descent method or Debye's method and the saddle-point
method \cite{copson:1965}. By rewriting $I_n$ in a suitable form,
we find that the saddle-point method provides the solution to our
asymptotic analysis.

Let
\begin{equation}\label{sec3-eq2}
g(t)=\mu(t)e^{t},
\end{equation}

then by the change of variables $t=nz$, our integral becomes
\begin{align}
I(n)&=n\int_{0}^{\infty}g(nz)
e^{-nz}\left\{\log\left(\frac{nz}{2\pi}\right)-i\frac{\pi}{2}\right\}^{n}\,dz\nonumber\\
&=n\int_{0}^{\infty}g(nz) e^{n\left\{-z+
\log\left[\log\left(\frac{nz}{2\pi}\right)-i\frac{\pi}{2}\right]\right\}}\,dz.\label{sec3-eq3}
\end{align}

If we define

\begin{equation}\label{sec3-eq4}
f(z)=-z+
\log\left[\log\left(\frac{nz}{2\pi}\right)-i\frac{\pi}{2}\right],
\end{equation}

then the saddle-point method consists in  deforming the path of
integration into a path which goes through a saddle-point at which
the derivative $f^{\prime}(z)$, vanishes. If $z_0$ is the
saddle-point at which the real part of $f(z)$ takes the greatest
value, the neighborhood of $z_0$ provides the dominant part of the
integral as $n\to \infty$ \cite[p. 91-93]{copson:1965}. This
dominant part provides an approximation of the integral and it is
given by the formula

\begin{align}
I(n)&\approx
ng(nz_0)e^{nf(z_0)}\left(\frac{-2\pi}{nf^{\prime\prime}(z_0)}
\right)^{\frac{1}{2}}.\label{sec3-eq5}
\end{align}

In our case, we have

\begin{align}
f^{\prime}(z)&=-1+\frac{1}{z\left[\log\left(\frac{nz}{2\pi}\right)-i\frac{\pi}{2}\right]},\quad\text{and}
\label{sec3-eq6}\\
f^{\prime\prime}(z)&=\frac{-1}{z^2\left[\log\left(\frac{nz}{2\pi}\right)-i\frac{\pi}{2}\right]}-
\frac{1}{z^2\left[\log\left(\frac{nz}{2\pi}\right)-i\frac{\pi}{2}\right]^2}.\label{sec3-eq7}
\end{align}

The saddle-point $z_0$ should verify the equation

\begin{align}
&\quad z_0\left[\log\left(\frac{nz_0}{2\pi}\right)-i\frac{\pi}{2}\right]=1\nonumber\\
\Leftrightarrow&\quad\frac{nz_0}{2\pi}\left[\log\left(\frac{nz_0}{2\pi}\right)-i\frac{\pi}{2}\right]=\frac{n}{2\pi}\nonumber\\
\Leftrightarrow&\quad\frac{nz_0}{2\pi}\log\left(\frac{nz_0}{2\pi}e^{-i\frac{\pi}{2}}\right)=\frac{n}{2\pi}\nonumber\\
\Leftrightarrow&\quad\frac{nz_0}{2\pi}e^{-i\frac{\pi}{2}}\log\left(\frac{nz_0}{2\pi}e^{-i\frac{\pi}{2}}\right)=
\frac{n}{2\pi}e^{-i\frac{\pi}{2}}.\label{sec3-eq8}
\end{align}

The last equation is of the form $v\log{v}=b$ whose solution can
be explicitly written using the principal branch\footnote{The
principal branch of the Lambert $W$-function is denoted by
$W_0(z)=W(z)$. See \cite{corless:lambert} for a thorough
explanation of the definition of all the branches.} of the Lambert
$W$-function \cite{corless:lambert}:

\begin{align}
v=e^{W(b)}.\label{sec3-eq9}
\end{align}

After some algebra, the saddle-point solution to our equation
(\ref{sec3-eq8}) is thus given by

\begin{align}
z_0=\frac{2\pi}{ni}e^{W\left(\frac{ni}{2\pi}\right)},\label{sec3-eq10}
\end{align}

and at the saddle-point, we have the values

\begin{align}\label{sec3-eq11}
f(z_0)&=-z_0-\log{z_0}\\
f^{\prime\prime}(z_0)&=-1-\frac{1}{z_0}.\label{sec3-eq12}
\end{align}

The saddle-point approximation of our integral (\ref{sec3-eq1}) is
given by the formula:

\begin{align}\label{sec3-eq13}
I(n)&=n\sqrt{\frac{2\pi}{n}}g(nz_0)e^{-nz_0-n\log(z_0)}\frac{1}{\sqrt{1+\frac{1}{z_0}}}\nonumber\\
&=n\sqrt{\frac{2\pi}{n}}\mu(nz_0)\frac{z_0^{\frac{1}{2}-n}}{\sqrt{1+z_0}}.
\end{align}

It turns out that $g(t)$ can be very well approximated by

\begin{equation}\label{sec3-eq14}
   g(t) =
\begin{cases} \frac{1}{6}e^{-\frac{1}{10}t^2} &\mbox{if } 0\le t\le 1 \\
-2+t & \mbox{if } t\gg 1.
\end{cases}
\end{equation}

Moreover, when $n$ is large, $g(nz_0)$ can also be very well
approximated\footnote{The approximations of $g(t)$ and $g(nz_0)$
are of course not necessary. We can keep the original functions
$g(nz_0)$ or $\mu(nz_0)$ for the final asymptotic formula.} by

\begin{equation}\label{sec3-eq15}
g(nz_0) \approx nz_0-1,
\end{equation}

so that we obtain the final approximation

\begin{align}\label{sec3-eq16}
I(n)\approx
n\sqrt{\frac{2\pi}{n}}(nz_0-1)\frac{z_0^{\frac{1}{2}-n}}{e^{nz_0}\sqrt{1+z_0}}.
\end{align}

To obtain an approximation of the coefficient
$\mu_n=\frac{2}{n!}\re\left\{I(n)\right\}$, we use Stirling
approximation of $n!$ and we further simplify $I(n)$  by resorting
to the following asymptotic development of the principal branch of
$W(z)$ \cite{corless:lambert}:

\begin{equation}\label{sec3-eq17}
W(z)=\log(z)-\log\left(\log{z}\right)+\cdots
\end{equation}

For $n\gg 1$, we can rewrite (\ref{sec3-eq10}) as

\begin{align}\label{sec3-eq18}
z_0&\sim \frac{1}{\log\left(\frac{n}{{2\pi}}\right)} e^{-i
\arctan\left(\frac{\pi}{2\log{n}}\right)}\sim
\frac{1}{\log\left(\frac{n}{{2\pi}}\right)} e^{-i
\frac{\pi}{2\log{n}}}
\end{align}
so that

\begin{align}\label{sec3-eq19}
\frac{1}{z_0^{n-\frac{1}{2}}}&\sim
\log\left(\frac{n}{{2\pi}}\right)^{n-\frac{1}{2}} e^{-i
(n-\frac{1}{2})\frac{\pi}{2\log{n}}},
\end{align}

and

\begin{align}\label{sec3-eq20} e^{-nz_0}&\sim
e^{-\frac{n}{\log\left(\frac{n}{{2\pi}}\right)}e^{-i\frac{\pi}{2\log{n}}}}\sim
e^{-\frac{n}{\log\left(\frac{n}{{2\pi}}\right)}} .
\end{align}

Using Stirling formula

\begin{align}\label{sec3-eq21}&
n!\sim \sqrt{2\pi n}\frac{n^n}{e^n},
\end{align}

 we obtain for large $n$

\begin{align}\label{sec3-eq22}
|\mu_n|&\sim \frac{n\log{n}}{e^{n\log{n}}}.
\end{align}

%%%%%%%%%%%%%%%%%%%%%%%%%%%%%%%%%%%%%%%%%%%%%%%%%%%%%%%%%%%%%%%%%%%%%%%%%%%%%%%%%%%%%%%%%%%%%%%%%%%%%%%%

The last equation proves that
$\lim_{n\to\infty}\frac{|\mu_{n+3}|}{|\mu_{n+2}|}=0$, or
equivalently that  $|\mu_{n+2}|\gg |\mu_{n+3}|\gg \ldots$ for
large $n$. Hence, Theorem~\ref{sec2-thm3} also provides an
asymptotic expansion of $\gamma_n$, and we deduce the following
one-term asymptotic approximation of $\gamma_{n}\approx
n!(-1)^n\mu_{n+2}$:

\begin{theorem}\label{sec3-thm1}
Let
$z_0^{*}=\frac{2\pi}{(n+2)i}e^{W\left(\frac{(n+2)i}{2\pi}\right)}$,
where $W$ is the Lambert $W$-function. An approximate formula for
the Stieltjes constants for large $n$ is
\begin{equation}\label{sec3-eq23}
\gamma_{n}\approx
\frac{2}{(n+1)}\sqrt{\frac{2\pi}{n+2}}\re\left\{\left((n+2)z_0^{*}-1\right)\frac{{z_0^{*}}^{\frac{1}{2}-n-2}}{e^{(n+2)z_0^{*}}
\sqrt{1+z_0^{*}}}\right\}.
\end{equation}
\end{theorem}

We can also find an asymptotic formula of $\gamma_n$ as a function
of $n$ only by using approximations similar to equations
(\ref{sec3-eq17}-\ref{sec3-eq19}). For $n\gg 1$ we can write

\begin{align}\label{sec3-eq24}
z_0^{*}&\sim  \frac{1}{\log\left(\frac{n+2}{{2\pi}}\right)} e^{-i
\frac{\pi}{2\log(n+2)}},
\end{align}

\begin{align}\label{sec3-eq25}
\frac{1}{{z_0^{*}}^{n+\frac{1}{2}}}&\sim
\log\left(\frac{n+2}{{2\pi}}\right)^{n+\frac{1}{2}} e^{-i
(n+\frac{1}{2})\frac{\pi}{2\log(n+2)}},
\end{align}

and

\begin{align}\label{sec3-eq26} e^{-(n+2)z_0^{*}}&\sim
e^{-\frac{(n+2)}{\log\left(\frac{n+2}{{2\pi}}\right)}} ,
\end{align}

and after some easy algebraic manipulations, we obtain the
oscillations and the leading order of growth of the Stieltjes
constants:

\begin{align}\label{sec3-eq27}
\gamma_n&\sim
2\frac{\sqrt{2\pi}}{\sqrt{n+2}}e^{(n+\frac{1}{2})\log(\log\left(n+2\right)-\log\left(2\pi\right))-
\frac{(n+2)}{\log\left(\frac{n+2}{{2\pi}}\right)}}\cos\left(
(n+\frac{1}{2})\frac{\pi}{2\log(n+2)}\right).
\end{align}

Both the oscillations and the leading order of growth match the
results of \cite{knessl:coffey}.

We also note that several terms of (\ref{sec2-eq18}) can also be
added to the one-term approximation given by
Theorem~\ref{sec3-thm1}. This leads to the following multi-term
approximation of $\gamma_n$:

\begin{theorem}\label{sec3-thm2}
Let
$z_0^{*}=\frac{2\pi}{(n+2)i}e^{W\left(\frac{(n+2)i}{2\pi}\right)}$,
where $W$ is the Lambert $W$-function. An M-term approximate
formula for the Stieltjes constants for large $n$ is
\begin{equation}\label{sec3-eq28}
\gamma_{n}\approx
\sum_{k=0}^{M-1}\frac{2n!}{(n+1+k)!}\sqrt{\frac{2\pi}{n+2+k}}\re\left\{\left((n+2+k)z_0^{*}-1\right)
\frac{{z_0^{*}}^{\frac{1}{2}-n-2-k}}{e^{(n+2+k)z_0^{*}}
\sqrt{1+z_0^{*}}}\right\}.
\end{equation}
\end{theorem}

When $M=1$, Theorem~\ref{sec3-thm2} reduces to
Theorem~\ref{sec3-thm1}. We will see in the next section that a
three-term ($M=3$) approximation provide satisfactory results for
large and small values of $n$.

\section{Numerical Results}\label{sec4}

We implemented the  formula of Theorem~\ref{sec3-thm2} in
Maple\textsuperscript{\scriptsize{\texttrademark}}
\blfootnote{\textsuperscript{\scriptsize{\texttrademark}}Maple is
a trademark of Waterloo Maple Inc.}. For a given value of $n$, the
following procedure  computes the value of the $M$-term
approximation formula of the $n^{\text{\tiny th}}$ Stieltjes
constant $\gamma_n$:
\begin{algorithmic}

\STATE \texttt{gamman := proc (n, M)}

\STATE \texttt{\#Input n: the desired nth Stieltjes constant}

\STATE \texttt{\#Input M: the number of terms in asymptotic
formula}

\STATE \texttt{\#An example call: gamman(137,3)}

\STATE \texttt{local k, coef, w0, z0, f, fpp;}

\STATE \texttt{coef := 0;}

\STATE \texttt{for k from 0 to M-1 do}

\STATE \texttt{w0 := LambertW(((1/2)*I)*(n+2+k)/Pi):}

\STATE \texttt{z0 := -(2*I)*Pi*exp(w0)/(n+2+k):}

\STATE \texttt{f := -z0-ln(z0):}

\STATE \texttt{fpp := -1-1/z0:}

\STATE \texttt{coef := coef+Re(2*factorial(n)*((n+2+k).z0-1)}

\STATE\texttt{*sqrt(-2*Pi/((n+2+k)*fpp))*exp((n+2+k)*f)/factorial(n+1+k)):}

\STATE \texttt{end do:}

\STATE \texttt{evalf(coef);}

\STATE \texttt{end proc:}

\end{algorithmic}

The approximations (\ref{sec3-eq23}) and (\ref{sec3-eq28}) with
$M=3$ were examined and compared to the exact values for $n$ from
$2$ to $100000$ given in
\cite{johansson:hurwitz,johansson:stieltjes} and the values of the
asymptotic formula of Knessl and Coffey \cite{knessl:coffey}.

Table~\ref{sec4-table1} below displays the approximate value of
$\gamma_{n}$ using Theorem~\ref{sec3-thm1} and
Theorem~\ref{sec3-thm2} and $M=3$, the approximation using the
formula of \cite{knessl:coffey} and the exact known values for $n$
from $2$ to $20$. Table~\ref{sec4-table2} displays the approximate
value of $\gamma_{n}$ and the exact known values for some higher
values of $n$.

\begin{table}[!h]
\centerline{
\begin{tabular}{ccccc}
  \hline
  $n$ & Exact & Theorem~\ref{sec3-thm2} &Theorem~\ref{sec3-thm1} &Knessl-Coffey\\
   & $\gamma_n$ & Eq. (\ref{sec3-eq28}), $M=3$  & Eq.  (\ref{sec3-eq23}) & Formula \cite{knessl:coffey}\\
  \hline
2&-0.009690363192&-0.008382380783&-0.008909030193&$\bf -$\\
3&0.002053834420 &0.001621242634&0.001073584137&0.00190188\\
4&0.002325370065 &0.002185636219&0.002025456323&0.00231644\\
5&0.000793323817 &0.0007895679944&0.000825888315&0.000812965\\
6&-0.000238769345&-0.0002241100338&-0.000149933239&-0.000242081\\
7&-0.000527289567&-0.0005200052551&-0.000475920788&-0.000541476\\
8&-0.000352123353&-0.0003506762586&-0.000346534072&-0.00036176\\
9&-0.000034394774&-0.0000349578308&-0.000055274760&-0.000035070\\
10&0.000205332814&0.0002044473764&0.000179950900&0.000210539\\
11&0.000270184439 &0.0002693789419&0.000255402785&0.00027624\\
12&0.000167272912 &0.0001666692377&0.000168701645&0.000170507\\
13&-0.000027463806 &-0.0000277087054&-0.000012840713&-0.000028263\\
14&-0.000209209262 &-0.0002089871741&-0.000190127572&-0.000213064\\
15&-0.000283468655 &-0.0002828583838&-0.000270364310&-0.000288108\\
16&-0.000199696858&-0.0001989876591& -0.000200475577&-0.000202633\\
17&0.000026277037 &0.0000266966357& 0.00000969746&0.0000267683\\
18&0.000307368408&0.0003071961365&0.000280749078&0.000311543\\
19&0.000503605453 &0.0005027990007&0.000479486029&0.000509981\\
20&0.000466343561 &0.0004652039644&0.000460162247&0.000471981\\
\hline
\end{tabular}
}\vspace{0.5cm} \caption{First 20 Stieltjes constants $\gamma_n$
and their approximate values given by Theorem~\ref{sec3-thm2},
Theorem~\ref{sec3-thm1}, and
by the formula of Knessl-Coffey.}%
\label{sec4-table1} %
\end{table}

\begin{table}[!h]
\centerline{
\begin{tabular}{ccccc}
  \hline
  $n$ & Exact & Theorem~\ref{sec3-thm2} &Theorem~\ref{sec3-thm1} &Knessl-Coffey\\
   & $\gamma_n$ & Eq. (\ref{sec3-eq28}), $M=3$  & Eq.  (\ref{sec3-eq23})& Formula \cite{knessl:coffey}\\
  \hline
30&0.003557728&0.0035491&  0.003790&0.00359535\\
35&-0.02037304&-0.0203320& -0.022336&-0.0205982\\
40&0.248721559&0.2484162& 0.265889&0.251108\\
45&-5.07234458 &-5.0686103&-5.211491&-5.10969\\
50& 126.8236026 &126.7545688&127.121&127.549\\
100&$-4.253401.10^{17}$ &$-4.251316.10^{17}$&$-4.14170.10^{17}$&$-4.25941.10^{17}$\\
136&$4.226701.10^{30}$ &$4.226998.10^{30}$&$4.22698.10^{30}$&$4.22698.10^{30}$\\
\bf 137& $\bf -0.00079.10^{29}$&$\bf -0.03484.10^{29}$&$ \bf 1.79099.10^{29}$&$\bf 3.89874.10^{29}$\\
138&$-2.523130.10^{31}$&$-2.521344.10^{31}$&$-2.4420176.10^{31}$&$-2.52354.10^{31}$\\
150&$8.028853.10^{35}$&$8.031999.10^{35}$&$8.1242241.10^{35}$&$8.05143.10^{35}$\\
250&$3.059212.10^{79}$ &$3.058889.10^{79}$&$3.038525.10^{79}$&$3.06165.10^{79}$\\
300&$-5.55672.10^{102}$ &$-5.55436.10^{102}$&$-5.47283.10^{102}$&$-5.55679.10^{102}$\\
800&$4.91354.10^{369}$ &$4.91329.10^{369}$&$4.899488.10^{369}$&$4.91452.10^{369}$\\
1400&$-4.09728.10^{728}$ &$-4.09772.10^{728}$&$-4.10081.10^{728}$&$-4.09851.10^{728}$\\
\hline
\end{tabular}
}\vspace{0.5cm} \caption{Stieltjes constants $\gamma_n$ and their
approximate values given by Theorem~\ref{sec3-thm2},
Theorem~\ref{sec3-thm1} and by Knessl-Coffey formula for different
values of $n$.}%
\label{sec4-table2} %
\end{table}

\begin{table}[!h]
\centerline{
\begin{tabular}{cc}
  \hline
  $n$ & Relative error \\
  \hline
2&-13.5~\%\\
3&-21.06~\%\\
4&-6.01~\%\\
6&-6.14~\%\\
137&-56.41~\% \\
821& 7.95~\%\\
1090& 8.01~\%\\
7259& 9.12~\%\\
8815&5.35~\%\\
\hline
\end{tabular}
} \vspace{0.5cm} \caption{The values of $n$ ($2\le n\le 100000$)
for which the relative error of the Stieltjes constants $\gamma_n$
computed
using Theorem~\ref{sec3-thm2} with $M=3$ exceeds 5\%.}%
\label{sec4-table3} %
\end{table}

We can see that both  asymptotic  formulas given by
Theorem~\ref{sec3-thm2} with $M=3$ or by Theorem~\ref{sec3-thm1}
provide  good approximations of the exact Stieltjes constants
except at $n=137$ where the approximation of
Theorem~\ref{sec3-thm1} fails to give the correct sign of
$\gamma_{137}$. Curiously, the asymptotic formula of Knessl-Coffey
also fails  to give the correct sign of $\gamma_{137}$. It seems
that the two asymptotic formulas are unrelated to each
other\footnote{The approximation formula of \cite{knessl:coffey}
is given by $\gamma_n\approx -\int_{0}^{\infty}\frac{\sin(\pi
e^{t})}{\pi} t^{n-1}e^{-t}(n-t)\,dt$, whereas our approximation
formula is $\gamma_n\approx
-\frac{1}{(n+1)(n+2)}\int_{0}^{\infty}2\mu(t)\re
\big\{(\log{t}-a)^{n+2}\big\}\,dt$. The author unsuccessfully
tried to derive a relationship between the two formulas.}. Thus,
the point $n=137$ is inherently  a badly conditioned point for
both asymptotic formulas.  For instance, with a small perturbation
of $n=137$, formula (\ref{sec3-eq23}) gives the value
$0.001041695409.10^{29}$ for $n=137.017$, and the value
$-0.1059515438.10^{29}$ for $n=137.018$. This shows that the point
$n=137$ is numerically ill-conditioned. This ill-conditioning can
be explained by the fact that the saddle-point equation of
\cite[eq. (2.4)]{knessl:coffey} and the saddle-point equation
(\ref{sec3-eq2}) both involve the evaluation of
$W(\frac{ni}{2\pi})$.

We also observe that except for the the specific values of $n=2$
and $n=137$, the approximation error of  $\gamma_n$ using the
formula of Knessl-Coffey is  less than that of the single-term
formula given by Theorem~\ref{sec3-thm1}. However, at the expense
of adding two extra terms to the approximation, the formula of
Theorem~\ref{sec3-thm2} outperforms both formulas: the sign error
of the $n=137$ disappears and the approximation error is greatly
reduced.

Table~\ref{sec4-table3} displays the values of $n$ ($2\le n\le
100000$) for which the relative error of the Stieltjes constants
$\gamma_n$ computed using Theorem~\ref{sec3-thm2} with $M=3$
exceeds 5\%. The relative errors rarely exceed $5\%$. In fact, for
$n\ge 8816$, the errors are  all less than $1.6\%$ with two
exceptions at $n=71158$ and $n=84589$ where the errors are equal
to $-2.4\%$ and $4.5\%$ respectively. It appears that with this
accuracy of the approximation, the three-term asymptotic formula
will hopefully be robust and work for all values of $n$.

\section{Conclusion and Extensions}\label{sec5}

It is possible\footnote{The functional equation artifice used in
this paper to find the asymptotic formula for $\gamma_n(a)$  when
$a=1$ extends easily to $a=\frac{1}{2}$. For irrational $a$ or
other rationals
 the extension is not obvious.} that the analysis of this paper can be generalized
to find  an effective asymptotic formulas for the generalized
Euler constants $\gamma_n(a)$ defined as the coefficients of the
Laurent series of the Hurwitz zeta function $\zeta(s,a)$ at the
point $s=1$ or at any other point of the complex plane. Instead of
formula (\ref{sec2-eq3}) of Theorem~\ref{sec2-thm1}, we use the
formula from \cite{fekih:hurwitz}:
\begin{equation}\label{sec4-eq1}
(s-1)\zeta(s,a) =\frac{1}{\Gamma(s)}\int_{0}^{\infty}\psi(t)
e^{-(a-1)t}t^{s-1}\,dt,
\end{equation}

which is valid for all $s$ such that $\re(s)>0$ and all $0<
a\le1$, and where the real function $\psi(t)$ is defined by

\begin{equation}\label{sec4-eq2}
\psi(t)=\frac{te^{t}}{(e^{t}-1)^2}-\frac{1}{e^{t}-1}+\frac{(a-1)t}{e^{t}-1}.
\end{equation}

It would be interesting to compare the formulas with the results
and conjectures of Kreminski  who has done extensive computations
on the generalized Euler constants \cite{kreminski:stieltjes}.

\section*{Acknowledgement}
I am grateful to Prof. Mark W. Coffey for his useful comments and
for providing the values of the last column in
Tables~\ref{sec4-table1}-\ref{sec4-table2} using the formula of
\cite{knessl:coffey}.

%%%%%%%%%%%%%%%%%%%%%%%%%%%%%%


\providecommand{\bysame}{\leavevmode\hbox to3em{\hrulefill}\thinspace}
\providecommand{\MR}{\relax\ifhmode\unskip\space\fi MR }
\providecommand{\MRhref}[2]{%
  \href{http://www.ams.org/mathscinet-getitem?mr=#1}{#2}
}
\providecommand{\href}[2]{#2}
\begin{thebibliography}{}

\end{thebibliography}


\begin{thebibliography}{00}

\bibitem{adell:stieltjes}J. A. Adell, Asymptotic estimates for Stieltjes constants. A
probabilistic approach, Proc. R. Soc. Lond. Ser. A Math. Phys.
Eng. Sci. 467, pp. 954–963, (2011).


\bibitem{briggs:zeta} W. E. Briggs, Some constants associated with the Riemann
zeta-function, Mich. Math. J. 3, pp. 117-121 (1955).

\bibitem{briggs:chowla}W. E. Briggs and S. Chowla, The power series coefficients of {$\zeta(s)$}, Amer. Math.
Monthly, Vol. 62, pp. 323-325, (1955).

\bibitem{coffey:stieltjes} M. W. Coffey, New results on the Stieltjes constants: Asymptotic and exact
evaluation, J. Math. Anal. Appl. Vol. 317, pp. 603–612, (2006).


\bibitem{copson:1965}E. T. Copson, Asymptotic Expansions, Cambridge University Press, (1965).

\bibitem{corless:lambert}R. M. Corles et al., On the Lambert W Function, Advances in Computational Mathematics,
Vol. 5, No. 1, pp. 329-359, (1996).

\bibitem{fekih:hurwitz}L. Fekih-Ahmed, On the Hurwitz Zeta Function,
\url{http://hal.archives-ouvertes.fr/hal-00602532}, (2011).


\bibitem{johansson:hurwitz}F. Johansson, Rigorous high-precision computation of the Hurwitz
zeta function and its derivatives,
\url{http://arxiv.org/abs/1309.2877}, (2013).

\bibitem{johansson:stieltjes}F. Johansson, The Stieltjes constants $\gamma_0 \ldots \gamma_{100000}$ rounded to 20 decimal
digits, \url{http://fredrikj.net/math/stieltjes100k20d.txt}, Jan.
(2014).


\bibitem{knessl:coffey}C. Knessl, M.W. Coffey, An effective asymptotic formula for the {S}tieltjes constants,
Mathematics of Computation, Vol. 80, No. 273, pp. 379-386, (2011).

\bibitem{kreminski:stieltjes}R. Kreminski, Newton-Cotes Integration for Approximating Stieltjes (Generalized Euler)
Constants, Mathematics of Computation, Vol. 72, No. 243, pp.
1379-1397, (2003).


\bibitem{matsuoka:zeta}Y. Matsuoka, On the power series coefficients of the Riemann zeta
function, Tokyo J. Math. 12, pp. 49-58 (1989).

\bibitem{mitrovic:stieltjes}D. Mitrovi\'{c}, The signs of some constants associated with the
Riemann zeta function, Mich. Math. J. 9, pp. 395-397 (1962).



\bibitem{sumaia:stieltjes}S. Saad Eddin, Explicit upper bounds for the Stieltjes
constants, Journal of Number Theory, Vol. 133, No. 3, pp.
1027-1044, (2013).


\bibitem{zhang:williams}N.-Y. Zhang and K. S. Williams, Some results on the generalized
Stieltjes constants, Analysis 14, pp. 147-162 (1994).


\end{thebibliography}
\end{document}